\documentclass[12pt]{amsart}
\input{amssymb.sty}
\title{A Kurosh type theorem for type $\mathrm{II}_1$ factors}
\author{Narutaka OZAWA}
\address{Department of Mathematical Sciences,
University of Tokyo, Komaba, 153-8914\\ 
Department of Mathematics, UCLA, Los Angeles, CA90095}
\email{narutaka@ms.u-tokyo.ac.jp}
\date{January 1, 2004}
\thanks{The author was supported by the JSPS 
Postdoctoral Fellowships for Research Abroad}

\newtheorem{thm}{Theorem}[section]
\newtheorem{prop}[thm]{Proposition}
\newtheorem{cor}[thm]{Corllary}
\newtheorem{lem}[thm]{Lemma}
\theoremstyle{definition}
\newtheorem{defn}[thm]{Definition}
\newtheorem{rem}[thm]{Remark}

\newcommand{\C}{{\mathbb C}} 
\newcommand{\B}{{\mathbb B}} 
\newcommand{\M}{{\mathbb M}} 
\newcommand{\N}{{\mathbb N}} 
\newcommand{\K}{{\mathbb K}} 
\newcommand{\F}{{\mathbb F}} 
\newcommand{\CC}{{\mathcal S}} 
\newcommand{\MM}{{\mathcal M}} 
\newcommand{\NN}{{\mathcal N}} 
 
\newcommand{\QQ}{{\mathcal Q}} 
\newcommand{\RR}{{\mathcal R}} 
\newcommand{\hh}{{\mathcal H}} 
\newcommand{\hk}{{\mathcal K}} 
\newcommand{\A}{{\mathcal A}} 
\newcommand{\U}{{\mathcal U}} 
 
\newcommand{\oo}{0}
\newcommand{\G}{\Gamma} 
\newcommand{\p}{\varphi}
\newcommand{\ce}{E} 
\newcommand{\red}{\mathrm{red}}
\newcommand{\Proj}{\mathop{\mathrm{Proj}}} 
\newcommand{\Prob}{\mathop{\mathrm{Prob}}} 
\newcommand{\Tr}{\mathop{\mathrm{Tr}}} 
\newcommand{\id}{\mathrm{id}} 
\newcommand{\conv}{\mathop{\overline{\mathrm{conv}}^w}}
\newcommand{\Aut}{\mathop{\mathrm{Aut}}}
\newcommand{\Ad}{\mathop{\mathrm{Ad}}\nolimits}
\newcommand{\supp}{\mathop{\mathrm{supp}}\nolimits}
\newcommand{\ip}[1]{\langle#1\rangle} 
\makeatletter
\def\freeprodsize@{\@setfontsize\freeprodsize\@xxpt\@xxpt}
\def\freeprod@{\mathop{\hbox{\freeprodsize@ $*$}}}
\def\freeprod{\freeprod@\displaylimits}
\makeatother
\newcommand{\bigast}{\freeprod}
\begin{document}
\begin{abstract}
We prove a Kurosh type theorem for free-product type $\mathrm{II}_1$ factors. 
In particular, if $M =L\F_2\bar{\otimes}\RR$, then the free-product 
type $\mathrm{II}_1$ factors $M*...*M$ are all prime and pairwise non-isomorphic. 
We also study the case of crossed product type $\mathrm{II}_1$ factors. 
This paper is a continuation of our previous papers \cite{solid}\cite{prime}, 
where the structure of (tensor products of) word hyperbolic group 
type~$\mathrm{II}_1$ factors was studied.
\end{abstract}
\maketitle
\section{Introduction}
The classification of type~$\mathrm{II}_1$ factors (of discrete groups)
was initiated by Murray and von Neumann \cite{mvn} who 
distinguished the hyperfinite type~$\mathrm{II}_1$ factor $\RR$ 
from the group factor $L\F_r$ of the free group $\F_r$ 
on $r\geq2$ generators. 
Thirty years later, Connes \cite{connes} proved uniqueness of 
the injective type~$\mathrm{II}_1$ factor. 
Thus, the group factor $L\G$ of an ICC amenable group $\G$ 
is isomorphic to the hyperfinite type~$\mathrm{II}_1$ factor $\RR$.
On the other hand, \emph{the} isomorphism problem of free group factors 
remains open. 
To solve this problem, 
Voiculescu invented free probability theory, 
which led to a number of deep results on the structure 
of free group factors (cf.\ the survey paper \cite{voiculescu}). 
Apart from these results and results of Connes \cite{connescras} 
and Cowling and Haagerup \cite{ch}, 
the classification of type~$\mathrm{II}_1$ factors 
has been vague by and large. 
Recently, however, a breakthrough came when 
Popa \cite{popab}\cite{popam} found that unitary conjugacy results 
can be deduced from existence of finite-dimensional bimodules 
and obtained quite precise classification theorems 
for certain classes of type~$\mathrm{II}_1$ factors. 
On the other hand, a $C^*$-algebraic method \cite{solid} was 
proved to be useful in study of type~$\mathrm{II}_1$ factors. 
These methods in combination yielded 
some prime factorization results in \cite{prime}. 

This paper is a continuation of \cite{solid} and \cite{prime}, 
where the structure of (tensor products of) word hyperbolic group 
type~$\mathrm{II}_1$ factors was studied. 
In this paper, we will study the structure of free-products 
and crossed products of certain type~$\mathrm{II}_1$ factors. 
A crucial ingredient of the argument is a computation 
of the kernels of certain morphisms on $C^*$-algebras. 
The idea of exploiting a `boundary' to compute such kernels 
is due to Skandalis \cite{skandalis} and developed by 
Higson and Guentner \cite{hg}. 
We will take advantage of this idea.
We denote by $\CC$ the class of countable discrete groups $\G$ 
such that the left and right translation action of $\G\times\G$ 
on the Stone-{\v C}ech remainder $\partial^\beta\G=\beta\G\setminus\G$ 
is amenable (see Section \ref{sec:crp} for details). 
The class $\CC$ was suggested by Skandalis and 
it contains all subgroups of word hyperbolic groups 
and discrete subgroups of connected simple Lie groups of rank one 
\cite{hg}\cite{skandalis}. 
The class $\CC$ also contains a group 
with an infinite amenable normal subgroup (cf.\ Corollary \ref{cor:wre}). 
The main result of \cite{solid} was solidity of 
the group factor $L\G$ of a group $\G$ in $\CC$. 
The $q$-Gaussian von Neumann algebras 
(for certain values of $q$) are other examples of solid factors. 
Indeed, solidity for certain values of $q$ 
was proved by Shlyakhtenko \cite{shlyakhtenko} 
while factoriality for all values of $q$ 
was proved by Ricard \cite{ricard}. 
The main results of \cite{prime} were 
unique prime factorization and rigidity of their tensor products. 

The main result of this paper 
is a Kurosh type theorem for a free-product 
of certain type~$\mathrm{II}_1$ factors. 
Although the theorem is not as precise as 
the original Kurosh theorem in group theory, 
it implies, for instance, that the iterated free-product 
type~$\mathrm{II}_1$ factors 
\[ 
L\F_\infty\ast(L\F_\infty\bar{\otimes}\RR)^{\ast n},\ n=1,2,\ldots
\]
are mutually non-isomorphic. 
This is a contrast to Dykema's theorem (Theorem 3.5 in \cite{dyk}) 
that $L\F_\infty\ast(L\F_\infty\bar{\otimes}L^\infty[0,1])^{\ast n}$
are all isomorphic. 
There is an obvious similarity between these results 
and the isomorphism problem of free group factors. 
In fact, according to Dykema and R{\u a}dulescu \cite{dr2}, 
isomorphism of free group factors would imply that 
$\MM_1\ast\MM_2=\MM_1\ast\MM_2\ast L\F_\infty$ 
for any type~$\mathrm{II}_1$ factors $\MM_1$ and $\MM_2$. 
The proof of the above theorem consists 
of an adaptation for free-product of 
the method developed in \cite{solid} and \cite{prime} 
(which we will review in Section \ref{sec:pre}) 
and Popa's work \cite{popao} on normalizers in a free-product. 

\begin{defn}
A type~$\mathrm{II}_1$ factor $\MM$ is \emph{prime} 
if $\MM\neq\MM_1\bar{\otimes}\MM_2$ 
for any type~$\mathrm{II}_1$ factors $\MM_1$ and $\MM_2$. 
A (finite) von Neumann algebra $\MM$ is \emph{solid} if 
for any diffuse von Neumann subalgebra $\A$, 
the relative commutant $\A'\cap\MM$ is injective. 
A (finite) von Neumann algebra is \emph{semisolid} if 
for any type~$\mathrm{II}_1$ von Neumann subalgebra $\QQ$, 
the relative commutant $\QQ'\cap\MM$ is injective. 
A von Neumann algebra is \emph{semiexact} if 
it contains a ultraweakly dense exact $C^*$-algebra.
\end{defn}

There are obvious implications; 
$\mbox{solid}\Rightarrow\mbox{semisolid}\Rightarrow\mbox{prime}$
for a non-injective type~$\mathrm{II}_1$ factor. 
We will see none of these implications is reversible. 
For a technical reason, the results in this paper are valid 
only for semiexact von Neumann algebras.
There are plenty of semiexact von Neumann algebras. 
A discrete group $\G$ is exact if and only if its reduced group $C^*$-algebra 
$C^*_\lambda \G$ is an exact $C^*$-algebra (cf.\ \cite{kw}). 
Thus, by definition, the group von Neumann algebra $L\G$ of 
a discrete exact group $\G$ is semiexact. 
The above mentioned Kurosh type theorem implies that 
every free-product of semiexact type~$\mathrm{II}_1$ factors is prime. 
This gives an example of prime type~$\mathrm{II}_1$ factors 
which are not semisolid. 

We also deal with crossed product and prove that 
the group-measure space von Neumann algebra $\G\ltimes L^\infty[0,1]$ 
of a measure-preserving action of a group $\G$ on 
the standard probability space $[0,1]$ is 
semisolid provided that the group $\G$ is in $\CC$. 
This generalizes Adams' theorem \cite{adams} that 
a measurable orbit equivalence relation of 
a non-amenable hyperbolic group is indecomposable. 
This also gives an example of semisolid type~$\mathrm{II}_1$ factors 
with the property $(\G)$. 
Note that a type~$\mathrm{II}_1$ factor with the property $(\G)$ 
cannot be solid by Proposition~7 in \cite{solid}. 

\medskip

\noindent\textbf{Acknowledgment.}
I would like to thank Nate Brown, Ken Dykema, Dimitri Shlyakhtenko 
and Georges Skandalis for many helpful discussions. 
In particular, I would like to thank Sorin Popa 
from whom I have learned a lot through the collaboration \cite{prime}. 
This research was carried out while I was visiting 
the University of California at Berkeley and 
the Penn State University under the support of 
the Japanese Society for the Promotion of Science 
Postdoctoral Fellowships for Research Abroad. 
\section{Conventions and Preliminary Background}\label{sec:pre}
For a discrete group $\G$, we denote by $\lambda$ 
(resp.\ $\rho$) the left (resp.\ right) regular representation 
of $\G$ on $\ell_2\G$. 
The reduced group $C^*$-algebra $C^*_\lambda\G$ 
(resp.\ $C^*_\rho\G$) is the $C^*$-subalgebra in $\B(\ell_2\G)$ 
which is generated by $\lambda(\G)$ (resp.\ $\rho(\G)$) 
and the group von Neumann algebra $L\G$ 
is the von Neumann algebra generated by $\lambda(\G)$. 
We denote by $C^*_{\lambda,\rho}\G$ the $C^*$-subalgebra 
in $\B(\ell_2\G)$ which is generated 
by $C^*_\lambda\G$ and $C^*_\rho\G$. 
Given a finite von Neumann algebra $\MM$, 
we assume that there is a distinguished 
faithful normal trace $\tau$ on $\MM$ 
and the trace on its von Neumann subalgebra $\NN\subset\MM$ 
is the restriction of $\tau$ on $\NN$. 
So, we will write $L^2\MM$ without specifying the trace $\tau$. 
We denote by $\widehat{a}$ the vector in $L^2\MM$ 
corresponding to $a\in\MM$. 
The Hilbert space $L^2\MM$ is an $\MM$-$\MM$ bimodule 
with $a\widehat{x}b=\widehat{axb}$ for $a,b,x\in\MM$. 
The canonical conjugation $J_{\MM}$ on $L^2\MM$ is given by 
$J\widehat{a}=\widehat{a^*}$ for $a\in\MM$. 
We simply denote $J_{\MM}$ by $J$ if there are no confusions. 
We will use the same notations for a $C^*$-algebra 
with a faithful trace. 
When dealing with $C^*$-algebras, 
the symbol $\otimes$ means the algebraic tensor product 
while $\otimes_{\min}$ means the minimal (spatial) tensor product. 
The spatial tensor product of von Neumann algebras is 
denoted by $\bar{\otimes}$ 
and the Hilbert space tensor product of Hilbert spaces 
is simply denoted by $\otimes$. 
The term `ucp' is an abbreviation for `unital completely positive'.
All representations and homomorphisms are assumed 
to be self-adjoint and non-degenerate. 
All groups denoted by $\G$ and $\Delta$ 
are assumed to be countable and discrete. 
All von Neumann algebras are assumed to have separable predual. 

We review the method developed in \cite{solid} and \cite{prime}. 
It tells when von Neumann subalgebra $\QQ$ 
in a finite von Neumann algebra $\MM$ has 
the relative commutant $\QQ'\cap\MM$ that is not injective. 

Let $\MM\subset\B(L^2\MM)$ be a finite von Neumann algebra 
and $\QQ\subset\MM$ be an injective von Neumann subalgebra. 
Then, there exists a conditional expectation 
$\Psi_{\QQ}$ from $\B(L^2\MM)$ onto $\QQ'$ which is proper, i.e.,
\[
\forall x\in\B(L^2\MM)\quad 
\Psi_{\QQ}(x)\in\conv\{ uxu^* : u\in\U(\QQ)\}.
\]
It follows that $\Psi_{\QQ}|_{\MM'}=\id_{\MM'}$ and 
that $\Psi_{\QQ}|_\MM$ is a trace preserving conditional expectation 
from $\MM$ onto $\QQ'\cap\MM$, which coincides with 
the unique trace preserving conditional expectation 
$\ce_{\QQ'\cap\MM}$ from $\MM$ onto $\QQ'\cap\MM$. 
Since $\Psi_{\QQ}|_{\MM'}=\id_{\MM'}$, the ucp map $\Psi_{\QQ}$ is 
an $\MM'$-bimodule map; 
\[
\forall x,y\in\MM',\ \forall b\in\B(L^2\MM)\quad 
\Psi_{\QQ}(xby)=x\Psi_{\QQ}(b)y.
\]
In particular, we have 
\[
\forall a\in\MM,\ \forall x\in\MM'\quad 
\Psi_{\QQ}(ax)=\ce_{\QQ'\cap\MM}(a)x.
\]
Now Lemma~5 in \cite{solid} can be interpreted as follows; 
\begin{lem}[\cite{solid}]\label{lem:inj}
Let $\MM\subset\B(L^2\MM)$ be a finite von Neumann algebra 
with an injective von Neumann subalgebra $\QQ\subset\MM$. 
If there exist unital ultraweakly dense $C^*$-subalgebras 
$B\subset\MM$ and $C\subset\MM'$ with $B$ exact such that 
the ucp map 
\[
\tilde{\Psi}_{\QQ}\colon B\otimes C\ni\sum_{k=1}^n a_k\otimes x_k
\mapsto\Psi_{\QQ}(\sum_{k=1}^n a_k x_k)\in\B(\hh)
\]
is continuous w.r.t.\ the minimal tensor norm on $B\otimes C$, 
then the relative commutant $\QQ'\cap\MM$ is injective.
\end{lem}

As in \cite{prime}, in an actual application, 
non-injectivity of $\QQ'\cap\MM$ forces 
\[
\K(\hk)\otimes_{\min}\B(L^2\NN)\not\subset\ker\Psi_{\QQ}
\] 
for some von Neumann subalgebra $\NN\subset\MM$ 
and a Hilbert space $\hk$ such that $L^2\MM=\hk\otimes L^2\NN$ 
as a right $\NN$ module.
If this is the case, we may find a finite rank projection 
$p$ on $\hk$ with $b=\Psi_{\QQ}(p\otimes 1_{\NN})\neq0$. 
Since $\Psi_{\QQ}$ is proper, 
$b$ commutes with the right $\NN$ action, 
or equivalently $b\in(\B(\hk)\bar{\otimes}\NN)\cap\QQ'$. 
By Proposition~1.3.2 in \cite{popab}, we have 
$(\Tr\otimes\tau_{\NN})(b) \leq \Tr(p)<\infty$.
Thus, there is a non-zero spectral projection $e$ of $b$ 
with $(\Tr\otimes\tau_{\NN})(e)<\infty$. 
It follows that $\hh=eL^2\MM$ is 
a $\QQ$-$\NN$ sub-bimodule of $L^2\MM$ with 
$\dim_{\NN}\hh_{\NN}<\infty$. 
(Strictly speaking, $\hh$ is 
a $\QQ e$-$\NN e'$ bimodule where $e'=J_{\MM}eJ_{\MM}$.)
If in addition $\MM$ and $\NN$ are factors, 
then we can apply the following 
Lemma~5 in \cite{popaf} and Proposition~12 in \cite{prime}. 

\begin{lem}[\cite{prime}\cite{popaf}]\label{lem:fdb}
Let $\NN$ and $\QQ$ be subfactors 
in a type~$\mathrm{II}_1$ factor $\MM$. 
If there exists a non-zero $\QQ$-$\NN$ sub-bimodule $\hh\subset L^2\MM$ 
with $\dim_{\NN} \hh_{\NN}<\infty$, then there exist 
projections $e\in\Proj(\NN)$ and $q\in\Proj(\QQ)$, 
a non-zero partial isometry $v\in\MM$ and 
a homomorphism $\theta\colon q\QQ q\to e\NN e$ such that 
\[
vv^*\in (q\QQ q)'\cap q\MM q,\ v^*v\in\theta(q \QQ q)'\cap e\MM e
\mbox{ and }xv=v\theta(x)\mbox{ for }x\in q\QQ q.
\]
\end{lem}

We need two more lemmas. 
The first is about normalizers in a free-product due to Popa. 
The following is formally stronger than stated in 
Theorem~6.1 and Corollary~4.3 in \cite{popao}, but 
their proofs are same (along Proposition~4.1 in \cite{popao}).
\begin{lem}[\cite{popao}]\label{lem:nif}
Let $\NN_1\subset\MM_1$ and $\MM_2$ be finite von Neumann algebras 
and let $\MM=\MM_1\ast\MM_2$ be their free-product. 
If $\NN_1$ is diffuse and a unitary operator $u\in\MM$ satisfies 
$u^*\NN_1u\subset\MM_i$ for some $i$, then $i=1$ and $u\in\MM_1$. 
\end{lem}

The last lemma in this section is 
about nuclearity of reduced free-product $C^*$-algebras. 
We recall that a ucp map $\p\colon A\to B$ is said to be 
\emph{nuclear} if there exist nets of 
ucp maps $\beta^\lambda\colon A\to\M_{n(\lambda)}(\C)$ 
and $\alpha^\lambda\colon\M_{n(\lambda)}(\C)\to B$ such that 
$\alpha^\lambda\circ\beta^\lambda\to\p$ in the point-norm topology. 
If $\p\colon A\to B$ is a nuclear ucp map and $B\subset\B(\hh)$, 
then the ucp map 
\[
\p\times\id_{B'}\colon A\otimes B'
\ni\sum a_k\otimes x_k\mapsto\sum\p(a_k)x_k\in\B(\hh)
\] 
is continuous w.r.t.\ the minimal tensor norm on $A\otimes B'$.

\begin{lem}\label{lem:nuc}
Let $B_i\subset\B(\hh_i)$ be a $C^*$-subalgebra with 
a $B_i$-cyclic unit vector $\xi_i\in\hh_i$.
We denote by $\omega_i$ the vector state corresponding to $\xi_i$.
If both $B_i$ are exact, then the inclusion map of
$(B_1,\omega_1)*(B_2,\omega_2)$ into $(\B(\hh_1),\omega_1)*(\B(\hh_2),\omega_2)$
is nuclear.
\end{lem}
\begin{proof}
It suffices to show that the inclusion map is approximated by 
ucp maps which factor through nuclear C$^*$-algebras.
We first note that if $B\subset\B(\hh)$ is exact,
then so is $B+\K(\hh)$.
(This is well-known and the proof is involved.
Indeed, by Kirchberg's theorem \cite{kirchberg}, 
$B$ is locally reflexive 
and the quotient $C:=B/(B\cap\K(\hh))$ is exact.
Moreover, the short exact sequence
\[
0\to\K(\hh)\to B+\K(\hh)\to C\to 0
\]
has ucp local splittings by the Effros-Haagerup lifting theorem \cite{eh}.
Now, the exactness of $B+\K(\hh)$ follows from that of $\K(\hh)$ and $C$
by the 3-by-3 lemma.)
By replacing $B_i$ with $B_i+\K(\hh_i)$, we may assume that 
$\K(\hh_i)\subset B_i$.
Since $B_i$ is exact, there are a net of finite dimensional subspaces
$\hk_i^\lambda\subset\hh_i$ with $\xi_i\in\hk_i^\lambda$
with the corresponding compression 
$\beta_i^\lambda\colon B_i\to\B(\hk_i^\lambda)$, and a net of ucp maps 
$\alpha_i^\lambda\colon\B(\hk_i^\lambda)\to\B(\hh_i)$
such that the net $\alpha_i^\lambda\circ\beta_i^\lambda$ converges
pointwise to the inclusion $B_i\hookrightarrow\B(\hh_i)$.
Since the rank-one projection $p_i$ corresponding to $\xi_i$
is in $B_i$, we have $\lim_\lambda \alpha_i^\lambda(p_i)=p_i$.
Thus, by perturbing $\alpha_i^\lambda$, we may assume that
$\alpha_i^\lambda(p_i)=p_i$ for all $\lambda$.
It follows from Choda-Blanchard-Dykema's theorem \cite{bd} that
\[
\beta_1^\lambda*\beta_2^\lambda\colon (B_1,\omega_1)*(B_2,\omega_2)
\to (\B(\hk_1^\lambda),\omega_1)*(\B(\hk_2^\lambda),\omega_2)
\]
and
\[
\alpha_1^\lambda*\alpha_2^\lambda\colon
(\B(\hk_1^\lambda),\omega_1)*(\B(\hk_2^\lambda),\omega_2)
\to (\B(\hh_1),\omega_1)*(\B(\hh_2),\omega_2)
\]
are ucp maps such that the net
$(\alpha_1^\lambda*\alpha_2^\lambda)\circ(\beta_1^\lambda*\beta_2^\lambda)$
converges pointwise to the inclusion map 
$(B_1,\omega_1)*(B_2,\omega_2)\hookrightarrow(\B(\hh_1),\omega_1)*(\B(\hh_2),\omega_2)$.
Since the C$^*$-algebras 
$(\B(\hk_1^\lambda),\omega_1)*(\B(\hk_2^\lambda),\omega_2)$
are nuclear (see e.g., \cite{kirchbergcom} or \cite{zzz}), we are done.
\end{proof}
\section{Free-Product}\label{sec:frp}
We recall the reduced free-product construction 
(in the tracial setting). 
Let $B_i$, $i\in\{1,2\}$ be a $C^*$-algebra 
with a faithful trace $\tau_i$. 
Let $B_i$ act on the GNS-Hilbert space $\hh_i=L^2(B_i,\tau_i)$. 
We denote by $\widehat{a}\in \hh_i$ the vector associated with 
$a\in B_i$ and denote by $\xi_i=\widehat{1}\in \hh_i$ 
the cyclic separating trace vector for $B_i$. 
The canonical conjugation $J_i$ on $\hh_i$ is given by 
$J_i\widehat{a}=\widehat{a^*}$.
Let $B_i^\oo=\ker\tau_i\subset B_i$ and 
let $\hh_i^\oo=\hh_i\ominus\C\xi_i$ be the closure of $B_i^\oo$ in $\hh_i$.
Then the free-product Hilbert space is 
\[
\hh=\C\xi\oplus\bigoplus_{n\geq1}\bigoplus_{\begin{subarray}{c}
i_1\ldots,i_n\in\{1,2\},\\ 
i_1\neq i_2,\cdots,i_{n-1}\neq i_n\end{subarray}}
\hh_{i_1}^\oo\otimes \hh_{i_2}^\oo\otimes\cdots\otimes \hh_{i_n}^\oo.
\]
We shall describe the left action $\lambda_i\colon\B(\hh_i)\to\B(\hh)$. 
It is convenient to introduce a subspace 
$\hh(i)\subset\hh$ which is the closed span of 
$\C\xi$ and those direct summands with $i\neq i_1$ 
in the above representation. 
Then, there is a canonical unitary operator 
\[
U_i\colon \hh \to \hh_i\otimes \hh(i)
\]
which identify $\hh(i)\cong\C\xi_i\otimes \hh(i)$ and 
$\hh(i)^\perp\cong\hh_i^\oo\otimes \hh(i)$.
We define 
\[
\lambda_i\colon\B(\hh_i)\ni a\mapsto 
U_i^* (a\otimes 1_{\hh(i)}) U_i \in\B(\hh).
\]
The reduced free-product $C^*$-algebra 
$(B,\tau)=(B_1,\tau_1)\ast(B_2,\tau_2)$ 
is the $C^*$-algebra $B$ in $\B(\hh)$ 
generated by $\lambda_1(B_1)$ and $\lambda_2(B_2)$ 
with the distinguished trace 
$\tau(\,\cdot\,)=(\,\cdot\,\xi,\xi)$ on $B$. 
We will omit $\lambda_i$ when there are no confusions.
The vector $\xi\in\hh$ is a cyclic separating trace vector for $B$ and 
the corresponding conjugation operator $J$ is given by 
\[
J(\widehat{a_1}\otimes\cdots\otimes\widehat{a_n})
=\widehat{a_n^*}\otimes\cdots\otimes\widehat{a_1^*}
=(J_{i_n}\widehat{a_n})\otimes\cdots\otimes(J_{i_1}\widehat{a_1})
\]
for $a_k\in B_{i_k}^\oo$ with $i_1\neq i_2,\cdots,i_{n-1}\neq i_n$.
In particular, $J\hh(i)\subset\hh$ is the closed linear span of $\C\xi$ 
and $\widehat{a_1}\otimes\cdots\otimes\widehat{a_n}$'s
with $i_n\neq i$. 
Let $V_i\colon\hh\to J\hh(i)\otimes\hh_i$ be the unitary operator 
which identifies $J\hh(i)\cong J\hh(i)\otimes\C\xi_i$
and $(J\hh(i))^\perp\cong J\hh(i)\otimes\hh_i^\oo$.
We note that $V_i$ intertwines the right actions of $\MM_i$. 
Moreover, the following is true. 

\begin{lem}\label{lem:mlt}
We have 
\begin{align*}
\lambda_i(a) &= JV_i^*(1_{J\hh(i)}\otimes J_iaJ_i)V_iJ\\
 &= V_i^*\bigl( P_\xi\otimes a 
  + \lambda_i(a)_{|J\hh(i)\ominus\C\xi}\otimes 1_{\hh_i}\bigr)V_i\\
 &= V_j^*(\lambda_i(a)_{|J\hh(j)}\otimes 1_{\hh_j})V_j
\end{align*}
for every $a\in\B(\hh_i)$ and $j\neq i$, where $P_\xi$ is 
the orthogonal projection onto $\C\xi$. 
\end{lem}
\begin{proof}
We first observe that the unitary operator 
\[
(J_{|J\hh(i)}\otimes J_i)V_iJU_i^*\colon\hh_i\otimes\hh(i)\to\hh(i)\otimes\hh_i
\]
is nothing but the flip of the tensor components $\hh_i$ and $\hh(i)$. 
Hence, we have $JV_i^*(J\otimes J_i)(1\otimes a)(J\otimes J_i)V_iJ
=U_i^*(a\otimes 1)U_i=\lambda_i(x)$ and the first equation follows. 
We denote by $P_i$ the orthogonal projection 
from $\hh$ onto $\C\xi\oplus\hh_i^\oo$.
We note that $P_i$ commutes with $\lambda_i(a)$. 
Since $U_iV_i^*(\xi\otimes\zeta)=\zeta\otimes\xi$ for every $\zeta\in\hh_i$, 
we have 
\[
V_i\lambda_i(a)P_iV_i^*(\xi\otimes\zeta)=V_iU_i^*(a\zeta\otimes\xi)
=\xi\otimes a\zeta=(P_\xi\otimes a)(\xi\otimes\zeta)
\]
for every $\zeta\in\hh_i$. 
Hence, we have $\lambda_i(a)P_i=V_i^*(P_\xi\otimes a)V_i$. 
We observe that 
\[
(U_i\otimes 1_{\hh_i})V_i(1-P_i)=(1_{\hh_i}\otimes V_i)U_i(1-P_i)
\]
as a partial isometry from $(1-P_i)\hh$ onto 
$\hh_i\otimes(\hh(i)\cap J\hh(i)\ominus\C\xi)\otimes\hh_i$. 
It follows that 
\begin{align*}
V_i^*(\lambda_i(a)\otimes 1)V_i(1-P_i) 
&= V_i^*(U_i^*\otimes 1)(a\otimes 1\otimes 1)(U_i\otimes 1)V_i(1-P_i)\\
&= V_i^*(U_i^*\otimes 1)(a\otimes 1\otimes 1)(1\otimes V_i)U_i(1-P_i)\\
&= V_i^*(U_i^*\otimes 1)(1\otimes V_i)U_i(1-P_i)\lambda_i(a)(1-P_i)\\
&= \lambda_i(a)(1-P_i).
\end{align*}
Since $Q_i:=V_i(1-P_i)V_i^*$ is the projection 
onto $(J\hh(i)\ominus\C\xi)\otimes\hh_i$, we have 
\[
\lambda_i(a)=\lambda_i(a)P_i+\lambda_i(a)(1-P_i)
=V_i^*\bigl((P_\xi\otimes a)+(\lambda_i(a)\otimes 1)Q_i\bigr)V_i
\]
as we claimed. 
Finally, if $i\neq j$, then $(U_i\otimes 1_{\hh_j})V_j=(1_{\hh_i}\otimes V_j)U_i$
and we have 
\[
V_j^*(\lambda_i(a)_{|J\hh(j)}\otimes 1)V_j
=U_i^*(1\otimes V_j^*)(a\otimes 1\otimes 1)(1\otimes V_j)U_i
=\lambda_i(a).
\]
\end{proof}
We denote $C_i=JB_iJ$, $C=JBJ$ and  
$D_i=V_i^*\bigl(\K(J\hh(i))\otimes_{\min}\B(\hh_i)\bigr)V_i$ 
for simplicity.
The above lemma implies that 
$\lambda_i(a)D_j\subset D_j$ for any $i,j\in\{1,2\}$ and 
$a\in\B(\hh_i)$.
\begin{prop}\label{prp:frp}
Let $(B,\tau)=(B_1,\tau_1)\ast(B_2,\tau_2)$, $\hh=L^2B$ and $C=JBJ$ 
be as above, and let $\Psi\colon\B(\hh)\to\B(\hh)$ 
be a $C$-bimodule ucp map. 
If $B_i$ are both exact and $D_i\subset\ker\Psi$ for both $i\in\{1,2\}$, 
then the ucp map 
\[
\tilde{\Psi}\colon B\otimes C\ni \sum_{k=1}^n a_k\otimes x_k
\mapsto \Psi(\sum_{k=1}^n a_k x_k) \in \B(\hh)
\]
is continuous w.r.t.\ the minimal tensor norm on $B\otimes C$.
\end{prop}
\begin{proof}
For simplicity, we let $\tilde{B}_i=\lambda(\B(\hh_i))$ and 
$\tilde{B}=C^*(\tilde{B}_1,\tilde{B}_2)$. 
We claim that $[\tilde{B},C]\subset\ker\Psi$. 
Since $\Psi$ is a $C$-bimodule map and 
the closed linear span of $\bigcup_{j=1,2}C[\tilde{B},C_j]C$ 
contains $[\tilde{B},C]$, it suffices to show 
$[\tilde{B},C_j]\subset\ker\Psi$ for each $j\in\{1,2\}$. 
But since $\tilde{B}D_k\tilde{B}\subset D_k$ by Lemma~\ref{lem:mlt} 
and the closed linear span of 
$\bigcup_{i=1,2}\tilde{B}[\tilde{B}_i,C_j]\tilde{B}$ 
contains $[\tilde{B},C_j]$, it suffices to show 
$[\tilde{B}_i,C_j]\subset\bigcup_{k=1,2} D_k$ 
for each $i,j\in\{1,2\}$. 
Now, it is not hard to see from Lemma~\ref{lem:mlt} that
$[\tilde{B}_i,C_j]=\{0\}$ for $i\neq j$ and that
\[
[\lambda_i(a),J\lambda_i(b)J] 
=V_i^*(\C P_\xi\otimes[a,J_ibJ_i])V_i\in D_i
\]
for every $a\in\B(\hh_i)$ and $b\in B_i$.

Since $[\tilde{B},C]\subset\ker\Psi$ and 
$\Psi$ is a $C$-bimodule map, we have $\Psi(\tilde{B})\subset C'$. 
Lemma~\ref{lem:nuc} implies that 
the inclusion map 
$\iota\colon B\hookrightarrow\tilde{B}$ 
is nuclear and so is the ucp map 
$\Psi_l=\Psi\circ\iota\colon B\to C'$. 
Therefore, the product ucp map 
\[
\tilde{\Psi}=\Psi_l\times\id_C\colon B\otimes C\to\B(\hh)
\]
is continuous w.r.t.\ the minimal tensor norm. 
\end{proof}

\begin{thm}\label{thm:frp}
Let $\MM_1$ and $\MM_2$ be semiexact finite factors 
and let $\MM=\MM_1\ast\MM_2$ be their free-product. 
If $\QQ\subset\MM$ is an injective type~$\mathrm{II}_1$ subfactor 
whose relative commutant $\QQ'\cap\MM$ is a non-injective factor, 
then there exist $i\in\{1,2\}$ and a unitary operator $u\in\MM$ 
such that $u^* \QQ u\subset \MM_i$ in $\MM$.
\end{thm}
\begin{proof}
We follow the notations used above. 
Let $B_i\subset\MM_i$ be ultraweakly dense exact $C^*$-algebras, 
$(B,\tau)=(B_1,\tau_1)\ast(B_2,\tau_2)$ be their free-product and 
let $C=JBJ$. 
Then, $B$ is ultraweakly dense in $\MM$ 
and is exact by Dykema's theorem \cite{dykema}.
It follows Lemma~\ref{lem:inj} that the $C$-bimodule ucp map 
\[
\tilde{\Psi}_{\QQ}\colon B\otimes C \ni \sum_{k=1}^n a_k\otimes x_k 
\mapsto \Psi_{\QQ}(\sum_{k=1}^n a_k x_k) \in \B(L^2\MM)
\]
cannot be continuous w.r.t.\ the minimal tensor norm. 
But by Proposition~\ref{prp:frp}, this implies that 
$D_i\not\subset\ker\Psi_{\QQ}$ for some $i\in\{1,2\}$. 
Now the discussion following Lemma~\ref{lem:inj} applies 
(for $\NN=\MM_i$) and yields a 
$\QQ$-$\MM_i$ sub-bimodule $\hk$ in $L^2\MM$ 
with $\dim \hk_{\MM_i}<\infty$. 

It follows Lemma~\ref{lem:fdb} that there exist 
projections $e\in\Proj(\MM_i)$ and $q\in\Proj(\QQ)$, 
a non-zero partial isometry $v\in\MM$ and 
a homomorphism $\theta\colon q\QQ q\to e\MM_i e$ such that 
$vv^*\in (q\QQ q)'\cap q\MM q$, $v^*v\in\theta(q \QQ q)'\cap e\MM e$ 
and $xv=v\theta(x)$ for $x\in q\QQ q$.
Since $(q\QQ q)'\cap q\MM q=q(\QQ'\cap\MM)q$, 
we have $vv^*=qq'$ for some $q'\in\Proj(\QQ'\cap\MM)$.
By restricting $\theta$ and $v$ if necessary, 
we may assume $\tau(q)=1/m$ and $\tau(q')=1/n$ 
for some $m,\,n\in\N$. 
Let $u_1,\ldots,u_m\in\QQ$ 
(resp.\ $u_1',\ldots,u_n'\in\QQ'\cap\MM$)
be partial isometries such that
$u_j^*u_j=q$ and $\sum_{j=1}^m u_j u_j^* =1$ 
(resp.\ $(u_k')^*u_k'=q'$ and $\sum_{k=1}^n u_k' (u_k')^* =1$).  
By Lemma~\ref{lem:nif}, $v^*v\in e\MM_ie$. 
We note that $\tau(v^*v)=\tau(qq')=(mn)^{-1}$. 
Let $w_{j,k}$ be partial isometries in $\MM_i$ such that 
$w_{j,k}w_{j,k}^*=v^*v$ and 
$\sum_{j=1}^m\sum_{k=1}^n w_{j,k}^*w_{j,k}=1$.
Then $u=\sum_{j,k}u_ju_k'vw_{j,k}$ is the desired unitary operator.
Indeed, 
$u^*xu
=\sum_{j_1,j_2,k}w_{j_1,k}^*\theta(u_{j_1}^*xu_{j_2})w_{j_2,k}\in\MM_i$ 
for $x\in\QQ$. 
\end{proof}

The Kurosh subgroup theorem states that 
if $\Lambda$ is a subgroup in a free product $\G_1\ast\G_2$, 
then $\Lambda$ is freely generated by a free subgroup in $\G_1\ast\G_2$
and/or conjugates of subgroups in $\G_1$ and/or $\G_2$. 
In particular, if $\Lambda\leq\G_1\ast\G_2$ 
is a freely-indecomposable non-infinite-cyclic subgroup, 
then $\Lambda$ is conjugated to a subgroup in $\G_1$ or $\G_2$. 
The following corollary is an analogue of this 
for type~$\mathrm{II}_1$ factors. 

\begin{cor}\label{cor:frp}
Let $\MM_1$ and $\MM_2$ be semiexact finite factors 
and let $\MM=\MM_1\ast\MM_2$ be their free-product. 
If $\NN\subset\MM$ is a non-prime non-injective 
subfactor whose relative commutant $\NN'\cap\MM$ is a factor, 
then there exist $i\in\{1,2\}$ and a unitary operator $u\in\MM$ 
such that $u^* \NN u \subset \MM_i$ in $\MM$.
In particular, $\MM$ is prime unless one of $\MM_i$ is trivial or 
both $\MM_i$ are isomorphic to $\M_2(\C)$.
\end{cor}
\begin{proof}
Since $\NN$ is non-prime non-injective, 
there are $\mathrm{II}_1$-factors $\NN_1$ and $\NN_2$ 
with $\NN_2$ non-injective such that 
$\NN=\NN_1\bar{\otimes}\NN_2$. 
By Proposition~13 in \cite{prime}, 
there is an injective type~$\mathrm{II}_1$ subfactor $\QQ\subset\NN_1$ 
with $\QQ'\cap\MM=\NN_1'\cap\MM$. 
Since $\NN_2\subset\NN_1'\cap\MM$, the relative commutant $\NN_1'\cap\MM$ 
is non-injective and the center of $\NN_1'\cap\MM$ 
is contained in (the center of) $\NN'\cap\MM$. 
Hence the factoriality of $\NN'\cap\MM$ implies that 
$\QQ'\cap\MM$ is a non-injective factor. 
Now it follows from Theorem~\ref{thm:frp} that 
there exist $i\in\{1,2\}$ and a unitary operator $u\in\MM$
such that $u^*\QQ u \subset\MM_i$. 
If $v\in\U(u^*\NN_2u)$, then $v$ commutes with $u^*\QQ u$ 
and hence $v\in\MM_i$ by Lemma~\ref{lem:nif}. 
It follows that $u^*\NN_2u\subset\MM_i$. 
The same argument applies to $u^*\NN_2u$ instead of $u^*\QQ u$ 
and we have $u^*\NN_1u\subset\MM_i$. 
Consequently, we have $u^*\NN u\subset \MM_i$.
We note that $\MM$ is a non-injective factor 
unless one of $\MM_i$ is trivial or 
both $\MM_i$ are isomorphic to $\M_2(\C)$ \cite{dyk}.
Since $\MM$ is not unitarily conjugated into $\MM_i$, 
it cannot be non-prime. 
\end{proof}

One of the consequences of the Kurosh subgroup theorem 
is the isomorphism theorem that 
if $\G_0$ and $\Lambda_0$ are free groups and 
$\G_1,\ldots,\G_n$ and $\Lambda_1,\ldots,\Lambda_m$ 
are freely-indecomposable non-infinite-cyclic groups 
with $\G=\bigast_{i=0}^n\G_i=\bigast_{j=0}^m\Lambda_j$, 
then $\G_0=\Lambda_0$, $n=m$ and, modulo permutation of indices, 
$\G_i$ and $\Lambda_i$ are conjugated in $\G$ for every $i\geq1$. 
The following is an analogue of this for type~$\mathrm{II}_1$ factors. 

\begin{cor}
Let $\MM_0,\ldots,\MM_n$ and $\NN_0,\ldots,\NN_m$ 
be semiexact finite factors such that 
$\MM_0$ and $\NN_0$ are semisolid 
(possibly one-dimensional) and that 
$\MM_1,\ldots,\MM_n$ and $\NN_1,\ldots,\NN_m$ are 
non-prime non-injective.
If $\MM=\bigast_{i=0}^n\MM_i=\bigast_{j=0}^m\NN_j$, 
then $n=m$ and, modulo permutation of indices, 
$\MM_i$ and $\NN_i$ are unitarily conjugated 
in $\MM$ for every $i\geq 1$. 
\end{cor}
\begin{proof}
Without loss of generality, we assume $m\geq n$ 
(and we no longer need the assumption that $\NN_0$ is semisolid).
Let $\MM=\bigast_{i=0}^n\MM_i=\bigast_{j=0}^m\NN_j$. 
By Corollary~\ref{cor:frp}, there exist maps 
$\imath\colon\{1,\ldots,m\}\to\{1,\ldots,n\}$ and 
$\jmath\colon\{1,\ldots,n\}\to\{0,\ldots,m\}$
and unitaries $u_1,\ldots,u_m$, $v_1,\ldots,v_n$ in $\MM$ 
such that 
$u_j^*\NN_ju_j\subset\MM_{\imath(j)}$ for every $j\geq1$
and $v_i^*\MM_iv_i\subset\NN_{\jmath(i)}$ for every $i\geq1$. 
It follows that 
$v_{\imath(j)}^*u_j^*\NN_ju_jv_{\imath(j)}\subset\NN_{\jmath(\imath(j))}$ 
for every $j\geq1$. 
But by Lemma~\ref{lem:nif}, this implies that $\jmath(\imath(j))=j$ and 
$u_jv_{\imath(j)}\in\NN_j$.
Therefore, the map $\imath$ is a bijection 
and the above inclusion maps are all isomorphisms. 
Hence, we have $u_j^*\NN_ju_j=\MM_{\imath(j)}$ for every $j\geq1$.
\end{proof}

\section{Crossed-Product}\label{sec:crp}
We recall the crossed product construction. 
Let $A$ be a $C^*$-algebra with a $\G$-action 
$\alpha\colon\G\to\Aut(A)$. 
A covariant representation of $(\G,A)$ is a pair $(\sigma,\pi)$ 
of representations of $\G$ and $A$ respectively 
on a Hilbert space $\hh$ such that 
\[
\forall s\in\G,\ \forall a\in A\quad
\Ad_{\sigma(s)}(\pi(a))=\sigma(s)\pi(a)\sigma(s)^{-1}=\pi(\alpha_s(a)).
\]
The full crossed product $C^*$-algebra $C^*(\G,A)$ 
is the universal $C^*$-algebra 
generated by $\sigma(\G)$ and $\pi(A)$ 
under the covariance condition. 
To define the reduced crossed product, 
we fix a faithful representation $\pi\colon A\to\B(\hh)$.
Define a new representation 
$\tilde{\pi}\colon A\to\B(\ell_2\G\otimes\hh)$ by 
\[
\tilde{\pi}(a)(\delta_t\otimes\zeta)
=\delta_t\otimes\pi(\alpha_t^{-1}(a))\zeta
\quad\mbox{for }t\in\G\mbox{ and }\zeta\in\hh.
\]
It is convenient to introduce orthogonal projections 
$e(t)$ of $\ell_2\G$ onto $\C\delta_t$ so that 
$\tilde{\pi}(a)=\sum_{t\in\G}e(t)\otimes\pi(\alpha_t^{-1}(a))$, 
where the sum converges strongly. 
It is easily verified that $(\lambda\otimes 1,\tilde{\pi})$ 
is a covariant representation of $(\G,A)$. 
The reduced crossed product $C^*$-algebra $C^*_{\red}(\G,A)$ 
is the $C^*$-subalgebra in $\B(\ell_2\G\otimes\hh)$ 
generated by $(\lambda\otimes1)(\G)$ and $\tilde{\pi}(A)$. 
We note that $C^*_{\red}(\G,A)$ does not depends on the choice 
of a faithful representation $\pi$. 
By definition, $C^*_{\red}(\G,A)$ is canonically isomorphic to 
a quotient of $C^*(\G,A)$. 
If $A$ is an exact $C^*$-algebra and $\G$ is an exact group, 
then $C^*_{\red}(\G,A)$ is exact also (cf.\ \cite{kw}).

A compact $\G$-space is a compact topological space $X$ together with 
a continuous action of $\G$ on it. 
Recall that we say a compact $\G$-space $X$ is \emph{amenable} 
(or, the $\G$-action on $X$ is amenable) 
if there exists a sequence of continuous 
$\mu^n\colon X\to\Prob(\G)$ such that 
\[
\lim_{n\to\infty}\sup_{x\in X}\|s.\mu^n_x-\mu^n_{s.x}\|=0
\]
for every $s\in\G$, where 
$\Prob(\G)=\{ \mu\in\ell_1\G : \mu\geq0,\ \|\mu\|=1\}$ 
and $(s.\mu)(t)=\mu(s^{-1}t)$ for $\mu\in\Prob(\G)$ and $s\in\G$. 
We note that the Stone-\v Cech compactification $\beta\G$ 
with the left translation action of $\G$ is amenable 
iff $\G$ is exact \cite{delaroche}\cite{gk}\cite{exact}. 
We assume that the reader is familiar with basic facts 
on amenability for group actions. 
We refer \cite{delaroche}\cite{adr} for detail. 

Let $X$ be a compact $\G$-space. 
A unital $\G$-$C(X)$-$C^*$-algebra 
is a unital $C^*$-algebra $A$ such that 
\begin{enumerate}
\item
$A$ contains $C(X)$ in its center, 
\item
there is a $\G$-action $\alpha\colon\G\to\Aut(A)$,
\item
$(\alpha_s(f))(x)=f(s^{-1}.x)$ for every $s\in\G$, $f\in C(X)$ and $x\in X$. 
\end{enumerate}
If the compact $\G$-space $X$ is amenable, 
then we have $C^*_{\red}(\G,A)=C^*(\G,A)$ canonically
for every $\G$-$C(X)$-$C^*$-algebra $A$. 

We first prove a general result on amenability of a group action. 
Suppose that $\G$ acts on a set $K$. 
The $\G$-action extends to a continuous action on 
the Stone-\v Cech compactification $\beta K$ of $K$ 
and then restricts to a continuous action on the Stone-\v Cech 
remainder $\partial^\beta K=\beta K\setminus K$. 
We are interested when the compact $\G$-space $\partial^\beta K$ is amenable. 
It is well-known that the unitary representation 
$\sigma_K\colon\G\to\B(\ell_2K)$, given by $\sigma_K(s)\delta_x=\delta_{s.x}$, 
is weakly contained in the left regular representation $\lambda$ iff 
the isotropy subgroups are all amenable. 
\begin{prop}\label{prop:pa}
Let $\G$ be an exact group, $K$ be a countable set on which $\G$ acts 
and $\partial^\beta K$ be the Stone-\v Cech remainder of $K$. 
The following are equivalent. 
\begin{enumerate}
\item\label{pa1} 
The compact $\G$-space $\partial^\beta K$ is amenable. 
\item\label{pa2} 
There exists a map 
$\mu\colon K\to\Prob(\G)$ such that 
$\lim_{x\to\infty}\|s.\mu_x-\mu_{s.x}\|=0$ for every $s\in\G$.
\item\label{pa3} 
There exists a ucp map $\p\colon C^*_\lambda\G\to\B(\ell_2K)$
such that $\p(\lambda(s))-\sigma_K(s)\in\K(\ell_2K)$ for every $s\in\G$. 
\end{enumerate}
\end{prop}
\begin{proof}
\ref{pa1}$\Rightarrow$\ref{pa3}: 
We note that 
$C(\partial^\beta K)\cong\ell_\infty K/c_0 K\subset\B(\ell_2 K)/\K(\ell_2 K)$. 
If the compact $\G$-space $\partial^\beta K$ is amenable, 
then the $C^*$-algebra $C^*_{\red}(\G,C(X))$ is nuclear 
and the natural homomorphism 
\[
C^*_\lambda\G\subset C^*_{\red}(\G,C(X))\cong C^*(\G,C(X))\to\B(\ell_2 K)/\K(\ell_2 K)
\]
is continuous. 
Moreover, it has a ucp lifting $\p$ by the Choi-Effros lifting theorem. 

\ref{pa3}$\Rightarrow$\ref{pa2}: 
By a generalized Weyl-von Neumann theorem (Theorem II.5.3 in \cite{davidson}), 
there exists an isometry $V\colon\ell_2K\to\ell_2\G$ such that 
$V^*\lambda(s)V-\p(\lambda(s))\in\K(\ell_2\G)$ for every $s\in\G$. 
For every $x\in K$, we set $\mu_x=|V\delta_x|^2\in\Prob(\G)$. 
Then, we have 
\begin{align*}
\|s.\mu_x-\mu_{s.x}\|_1 &\le \|\lambda(s)V\delta_x+V\delta_{s.x}\|_2
 \|\lambda(s)V\delta_x-V\delta_{s.x}\|_2\\
&\le 2(2-2\Re\ip{V^*\lambda(s)V\delta_x,\delta_{s.x}})^{1/2}.
\end{align*}
Since $V^*\lambda(s)V-\sigma_K(s)\in\K(\ell_2\G)$ for every $s\in\G$, we are done. 

\ref{pa2}$\Rightarrow$\ref{pa1}: 
Since the state space $S$ of $\ell_\infty\G$ is compact, 
the map $\mu\colon K\to\Prob(\G)$ extends to a 
continuous map $\tilde{\mu}$ from $\beta K$ into $S$. 
By the condition~\ref{pa2}, the map $\tilde{\mu}$ is $\G$-equivariant 
on $\partial^\beta K$. 
Hence, the amenability of $\partial^\beta K$ follows from that of $S$. 
However, since the state space $S$ is amenable iff 
the underlaying space $\beta\G$ is amenable, 
the amenability of $S$ follows from the exactness of $\G$. 
\end{proof}

Consider the $\G\times\G$-action on $\G$ given by the left and right 
translations. 
Let $\CC$ be the class of countable discrete groups $\G$ such that 
the compact $\G\times\G$-space $\partial^\beta\G$ is amenable. 
The class $\CC$ contains all word hyperbolic groups (and more generally 
groups which are hyperbolic relative to a family of amenable subgroups) 
\cite{hg}\cite{relhyp}\cite{skandalis}. 
Moreover, every subgroups of a group in $\CC$ is again in $\CC$ 
and $\CC$ is closed under free-product (with finite amalgamation).

Although it is irrelevant to the rest of paper, 
we make the following observation. 
A group $\G$ is said to be inner-amenable if 
there exists a state $m$ on $\ell_\infty\G$ with $c_o\G\subset\ker m$ 
such that $m(\sigma_s(f))=m(f)$ for every $s\in\G$ and $f\in\ell_\infty\G$, 
where $\sigma_s(f)(t)=f(s^{-1}ts)$. 
In other words, $\G$ is inner-amenable if 
$\partial^\beta\G$ carries a probability measure which is invariant 
under the conjugation action of $\G$. 
(There is another definition of inner-amenability 
that requires only $m(\delta_e)=0$ instead of $c_o\G\subset\ker m$. 
But, they coincide if the group in consideration is ICC.)
If $\G\in\CC$, then the conjugation action of $\G$ on $\partial^\beta\G$ 
is amenable since it is the restriction of 
the $\G\times\G$-action to its diagonal subgroup $\G$. 
It follows that every inner-amenable group in $\CC$ 
is amenable. 
This was first proved in \cite{delaharpe}. 
It is also not difficult to prove that 
every torsion-free non-amenable group in $\CC$ is ICC. 

We look at the crossed product construction more carefully 
in the tracial setting. 
Let $A$ be a $C^*$-algebra with a faithful trace $\tau$ 
and let $\alpha\colon\G\to\Aut(A,\tau)$ be a trace preserving action. 
If $\pi_\tau\colon A\to\B(L^2A)$ is 
the GNS-representation of $(A,\tau)$,  
then $C^*_{\red}(\G,A)$ is the $C^*$-subalgebra 
in $\B(\ell_2\G\otimes L^2A)$ generated by 
$(\lambda\otimes1)(\G)$ and $\tilde{\pi}_\tau(A)$. 
The vector $\xi=\delta_e\otimes\widehat{1}_A$ 
is a cyclic separating trace vector for $C^*_{\red}(\G,A)$ 
with the canonical conjugation $J$ on $\ell_2\G\otimes L^2A$ given by 
$Jx\xi=x^*\xi$ for $x\in C^*_{\red}(\G,A)$.
Let $u$ be the representation of $\G$ on $L^2A$ given by 
$u(s)\widehat{a}=\widehat{\alpha_s(a)}$ for $s\in\G$ and $a\in A$.
A simple calculation shows that 
\[
J(\lambda\otimes 1)(s)J=(\rho\otimes u)(s)
\mbox{ for }s\in\G\mbox{ and }
J\tilde{\pi}_\tau(a)J=1\otimes\pi_\tau^c(a) 
\mbox{ for }a\in A,
\]
where $\pi_\tau^c(a)=J_A\pi_\tau(a)J_A$. 
For simplicity, denote by $C^*(B,C)$ 
the $C^*$-subalgebra in $\B(\ell_2\G\otimes L^2A)$ 
generated by $B=C^*_{\red}(\G,A)$ and $C=JBJ$. 

Let $I_1=\K(\ell_2\G)\otimes_{\min}\B(L^2A)$.
We note that $I_1$ is the hereditary $C^*$-subalgebra 
in $\B(\ell_2\G\otimes L^2A)$. 
It is not hard to see that both $B$ and $C$ are 
in the multiplier of $I_1$ and hence 
$I=I_1\cap C^*(B,C)$ is an ideal in $C^*(B,C)$. 
\begin{prop}\label{prp:crp}
Let $\alpha\colon \G\to\Aut(A,\tau)$ 
and $I\subset C^*(B,C)$ be as above. 
Suppose that $A''\subset\B(L^2A)$ is injective and $\G\in\CC$. 
Then, the homomorphism 
\[
\nu\colon B\otimes C \ni \sum_{k=1}^n a_k\otimes x_k 
\mapsto \sum_{k=1}^n a_k x_k + I \in C^*(B,C)/I
\]
is continuous w.r.t.\ the minimal tensor norm on $B\otimes C$.
\end{prop}
\begin{proof}
Consider the $C^*$-algebra $D_1=C^*(B,C,\ell_\infty\G\otimes\C1)$ 
generated by $C^*(B,C)$ and $\ell_\infty\G\otimes1$. 
Since $D_1$ is in the multiplier of $I_1$, 
there is a natural inclusion
\[
C^*(B,C)/I\hookrightarrow(D_1+I_1)/I_1=E_1.
\]
Let $D\subset D_1$ be the $C^*$-algebra 
generated by $\tilde{\pi}(A)$, $J\tilde{\pi}(A)J$ 
and $\ell_\infty\G\otimes\C1$. 
Then, it is not hard to see that $E=(D+I_1)/I_1\subset E_1$ is 
a $(\G\times\G)$-$\partial\G$-$C^*$-algebra 
with the commuting $\G$-actions 
$\Ad_{\lambda\otimes 1}$ and $\Ad_{\rho\otimes u}$. 
Since $\G\in\CC$, 
the canonical homomorphism from 
$C^*(\G\times\G,E)$ onto $E_1$ 
factors through $C^*_{\red}(\G\times\G,E)$. 
Since $\tilde{\pi}(A)''$ is injective, the natural homomorphism 
$\tilde{\pi}(A)\otimes J\tilde{\pi}(A)J\to E$ is 
continuous w.r.t.\ the minimal tensor norm. 
Therefore, the homomorphism $\nu$ is continuous on 
$B\otimes_{\min}C\cong 
C^*_{\red}(\G\times\G, \tilde{\pi}(A)\otimes_{\min}J\tilde{\pi}(A)J)$.
\end{proof}

Let 
\[
K_1=\{x\in\B(\ell_2\G\otimes L^2A) : 
(\omega\otimes\id)(x^*x+xx^*)\in\K(L^2A)\ 
\forall\omega\in\B(\ell_2\G)_*^+\}.
\] 
We note that $K_1$ is the hereditary $C^*$-subalgebra 
in $\B(\ell_2\G\otimes L^2A)$ 
generated by $\ell_\infty(\G,\K(L^2A))$. 
It is not hard to see that both $B$ and $C$ are 
in the multiplier of $K_1$ and hence 
$K=K_1\cap C^*(B,C)$ is an ideal in $C^*(B,C)$. 

\begin{prop}\label{prp:kas}
Let $\alpha\colon \G\to\Aut(A,\tau)$ and $K\subset C^*(B,C)$ 
be as above. 
Suppose that $A''\subset\B(L^2A)$ is injective and 
there exists a unital completely positive map 
$\theta\colon\B(L^2A\otimes \ell_2\G)\to\B(L^2A)$ 
such that the elements
$\theta((u\otimes\rho)(s))-u(s)$, 
$\theta(\pi_\tau(a)\otimes1)-\pi_\tau(a)$ and 
$\theta(\pi_\tau^c(a)\otimes1)-\pi_\tau^c(a)$ 
are all in $\K(L^2A)$ for every $s\in\G$ and $a\in A$.
Then, the homomorphism 
\[
\mu\colon B\otimes C \ni \sum_{k=1}^n a_k\otimes x_k 
\mapsto \sum_{k=1}^n a_k x_k + K \in C^*(B,C)/K
\]
is continuous w.r.t.\ the minimal tensor norm on $B\otimes C$.
\end{prop}

\begin{proof}
Since $A''\subset\B(L^2A)$ is injective, 
there is a canonical homomorphism 
\[
\Phi\colon 
(\B(\ell_2\G)\bar{\otimes}A'')
\otimes_{\min}(A'\bar{\otimes}\B(\ell_2\G))
\to\B(\ell_2\G\otimes L^2A\otimes\ell_2\G)
\]
given by 
$\Phi(x\otimes 1\otimes 1)=x\otimes 1$ 
for $x\in \B(\ell_2\G)\bar{\otimes}A''$ 
and $\Phi(1\otimes1\otimes y)=1\otimes y$ 
for $y\in A'\bar{\otimes}\B(\ell_2\G)$.
Let $U_0\colon \ell_2\G\otimes L^2A\to L^2A\otimes\ell_2\G$ 
be the unitary operator given by 
$U_0(\delta_t\otimes\xi)\mapsto u(t)\xi\otimes\delta_t$
so that the elements
\[
\Ad_{U_0}((\rho\otimes u)(s))=(1\otimes\rho)(s),\
\Ad_{U_0}(1\otimes\pi_\tau^c(a))
=\sum_{t\in\G}\pi_\tau^c(\alpha_t(a))\otimes e(t)
\]
are in $A'\bar{\otimes}\B(\ell_2\G)$. 
Let $U_1=\sum_{t\in\G}(\rho\otimes u)(t)\otimes e(t)
\in B'\bar{\otimes}\B(\ell_2\G)$ 
be the unitary operator on $\ell_2\G\otimes L^2A\otimes\ell_2\G$
so that 
\[
\Ad_{U_1^*}((1\otimes 1\otimes \rho)(s))=(\rho\otimes u\otimes\rho)(s),\
\Ad_{U_1^*(1\otimes U_0)}(1\otimes 1\otimes\pi_\tau^c(a))
=1\otimes\pi_\tau^c(a)\otimes 1.
\] 
It follows that for the homomorphism 
\[
\tilde{\Phi}=\Ad_{U_1^*}\Phi\Ad_{1\otimes1\otimes U_0}
\colon B\otimes_{\min}C
\to\B(\ell_2\G\otimes L^2A\otimes\ell_2\G),
\]
we have 
\begin{align*}
\tilde{\Phi}((\lambda\otimes 1\otimes 1\otimes 1)(s))
&= \lambda(s)\otimes1\otimes1,\\
\tilde{\Phi}(\tilde{\pi}_\tau(a)\otimes 1\otimes 1)
&=\tilde{\pi}_\tau(a)\otimes1,\\
\tilde{\Phi}((1\otimes 1\otimes \rho\otimes u)(s))
&=(\rho\otimes u\otimes\rho)(s),\\
\tilde{\Phi}(1\otimes1\otimes1\otimes\pi_\tau^c(a))
&=1\otimes\pi_\tau^c(a)\otimes 1.
\end{align*}
It follows that for the unital completely positive map 
\[
\tilde{\theta}=\id\otimes\theta\colon
\B(\ell_2\G\otimes L^2A\otimes\ell_2\G)
\to\B(\ell_2\G\otimes L^2A), 
\]
the elements
\begin{align*}
&\tilde{\theta}\tilde{\Phi}((\lambda\otimes 1\otimes 1\otimes 1)(s))
-\lambda(s)\otimes1
=0,\\
&\tilde{\theta}\tilde{\Phi}(\tilde{\pi}_\tau(a)\otimes 1\otimes 1)
-\tilde{\pi}_\tau(a)
=\sum_{t\in\G}e(t)\otimes(\theta(\pi_\tau(\alpha_t^{-1}(a))\otimes 1)
-\pi_\tau(\alpha_t^{-1}(a))),\\
&\tilde{\theta}\tilde{\Phi}((1\otimes 1\otimes \rho\otimes u)(s))
-(\rho\otimes u)(s)
=\rho(s)\otimes(\theta((u\otimes\rho)(s))-u(s)),\\
&\tilde{\theta}\tilde{\Phi}(1\otimes1\otimes1\otimes\pi_\tau^c(a))
-1\otimes\pi_\tau^c(a)
=1\otimes(\theta(\pi_\tau^c(a)\otimes 1)-\pi_\tau^c(a))
\end{align*}
are all in $K_1$ for every $s\in\G$ and $a\in A$. 
Since $(C^*(B,C)+K_1)/K_1=C^*(B,C)/K$ canonically, 
the completely positive map 
$\tilde{\theta}\tilde{\Phi}\colon B\otimes_{\min}C\to C^*(B,C)+K_1$
passes to the homomorphism $\mu$.
\end{proof}

Let us recall the Bernoulli shift (or the wreath product) 
of $\Delta$ by $\G$ is the group $\G\ltimes\Delta_\G$, 
where $\Delta_\G=\bigoplus_\G\Delta$ 
is the direct sum of $\Delta$'s indexed by $\G$ 
and $\G$ acts on $\Delta_\G$ by left translation. 
We view an element $x\in\Delta_\G$ as a function 
$x\colon \G\to\Delta$ with finite support 
$\{ t\in\G : x(s)\neq e_\Delta\}=:\supp x$. 

\begin{prop}\label{prp:ber}
Let $\G$ and $\Delta$ be groups with $\Delta$ amenable. 
Then, the Bernoulli shift action $\G$ on $A=C^*_\lambda\Delta_\G$ 
satisfies the assumptions in Proposition \ref{prp:kas}.
\end{prop}

\begin{proof}
We note that $L^2A$ is canonically isomorphic to $\ell_2\Delta_\G$. 
We fix a proper length function $l_\G$ on $\G$, i.e., 
$l_\G$ is a non-negative function on $\G$ 
such that (i) $l_\G(s)=0$ iff $s=e$, 
(ii) $l_\G(st)\le l_\G(s)+l_\G(t)$ for $s,t\in\G$, 
and (iii) the set $\{s\in\G : l_\G(s)\le R\}$ 
is finite for every $R>0$. 
Likewise $l_\Delta$. 
For $y\in\Delta_\G$ and $t\in\G$, 
we set $w(y,t)=l_\G(t)+l_\Delta(y(t))$ if $t\in\supp y$ 
and $w(y,t)=0$ otherwise. 
Further, set $w(y)=\sum_{t\in\G}w(y,t)$ and $n(y)=|\supp y|$. 
It follows that $n(y)/w(y)\to0$ as $y\to\infty$. 
Define $\xi\colon\Delta_\G\to\ell_2\G$ by 
$\xi_y(t)=(w(y,t)/w(y))^{1/2}$. 
Since $|w(y,t)-w(s(y),st)|\le l_\G(s)$ 
for $s\in\G$ and $t\in\supp y$, 
we have for every $s\in\G$ that 
\begin{align*}
\|\lambda(s)\xi_y-\xi_{s(y)}\|_2^2
&\le\|(\lambda(s)\xi_y)^2-\xi_{s(y)}^2\|_1\\
&=\sum_{t\in\G}|w(y,t)/w(y) - w(s(y),st)/w(s(y))|\\
&=\left(\sum_{t\in\G}|w(y,t)-w(s(y),st)|/w(y)\right)+|w(s(y))-w(y)|/w(y)\\
&\le 2l_\G(s)n(y)/w(y)\to 0\mbox{ as }y\to\infty.
\end{align*}
Moreover for each $x\in\Delta_\G$, 
we have 
$\|\xi_{xy}-\xi_y\|_2^2\le 2w(x)/w(y)\to 0$ 
and $\|\xi_{yx}-\xi_y\|_2^2\le 2w(x)/w(y)\to 0$
as $y\to\infty$.

Let $V\colon \ell_2\Delta_\G\to\ell_2\Delta_\G\otimes\ell_2\G$ 
be the isometry given by 
$V\delta_y=\delta_y\otimes\xi_y$.
Then the unital completely positive map 
$\theta\colon
\B(\ell_2\Delta_\G\otimes \ell_2\G)\to\B(\ell_2\Delta_\G)$, 
given by $\theta(z)=V^*zV$, 
satisfies the assumption of Proposition \ref{prp:kas} 
(with $\rho$ replaced with $\lambda$). 
Indeed, 
\begin{align*}
V^*(u\otimes\lambda)(s)V\delta_y
&=\ip{\lambda(s)\xi_y,\xi_{s(y)}}\delta_{s(y)}\\
V^*(\lambda(x)\otimes 1)V\delta_y
&=\ip{\xi_y,\xi_{xy}}\delta_{xy}\\
V^*(\rho(x)\otimes 1)V\delta_y
&=\ip{\xi_y,\xi_{yx^{-1}}}\delta_{yx^{-1}}
\end{align*}
for every $s\in\G$ and $x\in\Delta_\G$.
\end{proof}

The group $\G$ acts on $\Delta_\G$ by Bernoulli shift and 
$\Delta_\G\times\Delta_\G$ acts on $\Delta_\G$ 
by left and right translations. 
These actions induce an action of 
$\Lambda:=\G\ltimes(\Delta_\G\times\Delta_\G)$ 
on the Stone-\v Cech remainder $\partial^\beta\Delta_\G$. 
We observe from the above proof that this action is amenable 
provided that $\G$ is exact. 
Indeed, the map $\xi^2\colon \Delta_\G\to\ell_1\G$ 
gives rise to a continuous map from 
the Stone-\v Cech compactification $\beta\Delta_\G$ 
into the state space of $\ell_\infty\G$ 
whose restriction to $\partial^\beta\Delta_\G$ 
is $\Lambda$-equivariant (where $\Delta_\G\times\Delta_\G$ 
acts trivially on $\ell_\infty\G$). 
Since $\Delta_\G$ is amenable there exists 
a $\Lambda$-equivariant conditional expectation from 
$\ell_\infty\Lambda$ onto $\ell_\infty\G$. 
Composing these maps, we obtain a $\Lambda$-equivariant 
continuous map from $\partial^\beta\Delta_\G$ into 
the state space of $\ell_\infty\Lambda$. 
Therefore, the amenability of $\partial^\beta\Delta_\G$ follows 
from that of the state space of $\ell_\infty\Lambda$, which 
is amenable when $\G$ (and hence $\Lambda$) is exact. 

\begin{cor}\label{cor:wre}
If $\G\in\CC$ and $\Delta$ is amenable, 
then the Bernoulli shift $\G\ltimes\Delta_\G$ is in $\CC$. 
\end{cor}
\begin{proof}
By Propositions \ref{prp:crp}, \ref{prp:kas} and \ref{prp:ber} 
(and their proof), the homomorphism 
\[
B\otimes C \ni \sum_{k=1}^n a_k\otimes x_k 
\mapsto \sum_{k=1}^n a_k x_k + (I\cap K) \in C^*(B,C)/(I\cap K)
\]
is continuous w.r.t.\ the minimal norm and 
has a ucp lifting. 
Since 
$I_1\cap K_1=\K(\ell_2\G\otimes\ell_2\Delta_\G)$, 
the claim follows from Proposition \ref{prop:pa}. 
\end{proof}

If $A=\A$ is a von Neumann algebra, then 
the crossed product von Neumann algebra is 
$\G\ltimes\A=C^*_{\red}(\G,\A)''\subset\B(\ell_2\G\otimes L^2\A)$.
The following are the main results of this section. 
\begin{thm}\label{thm:crp}
Let $(\A,\tau)$ be a commutative von Neumann algebra with 
a faithful trace $\tau$ and a trace preserving action 
$\alpha\colon\G\to\Aut(\A,\tau)$ and let $\MM=\G\ltimes\A$ 
be its crossed product (which may not be a factor). 
If $\G\in\CC$, then $\MM$ is semisolid.
In particular, any non-injective subfactor of $\MM$ is prime. 
\end{thm}
\begin{proof}
We follow the notations used above. 
In particular, $B=C^*_{\red}(\G,A)$ and $C=JBJ$. 
Let $\QQ\subset\MM$ be a type~$\mathrm{II}_1$ von Neumann subalgebra. 
Passing to a subalgebra if necessary, we may assume $\QQ$ is injective. 
For a proof by contradiction, 
suppose that $\QQ'\cap\MM$ is not injective. 
It follows from Lemma~\ref{lem:inj} that 
the ucp map 
\[
\tilde{\Psi}_{\QQ}\colon B\otimes C \ni \sum_{k=1}^n a_k\otimes x_k 
\mapsto \Psi_{\QQ}(\sum_{k=1}^n a_k x_k) \in \B(L^2\MM)
\]
cannot be continuous w.r.t.\ the minimal tensor norm. 
But, by Proposition~\ref{prp:crp}, this implies that 
$\K(\ell_2\G)\otimes\B(L^2\A)\not\subset\Psi_{\QQ}$.
Now the discussion following Lemma~\ref{lem:inj} 
applies (for $\NN=\A$) and yields 
a non-zero $\QQ$-$\A$ bimodule $\hh$ with $\dim\hh_\A<\infty$. 
This is absurd since $\A$ is commutative while $\QQ$ 
is of type~$\mathrm{II}_1$.
\end{proof}

\begin{thm}
Let $\G$ be an exact group and 
$\RR$ be the hyperfinite type~$\mathrm{II}_1$ factor. 
Consider the Bernoulli product 
$\MM=\G\ltimes\bar{\bigotimes}_\G\RR$. 
Then, for any diffuse von Neumann subalgebra 
$\QQ\subset\bar{\bigotimes}_\G\RR$, the relative commutant 
$\QQ'\cap\MM$ is injective. 
\end{thm}
\begin{proof}
Let $\Delta$ be an amenable ICC group so that $L\Delta=\RR$ 
and let $A=C^*_\lambda\Delta_\G$. 
Let a diffuse von Neumann subalgebra 
$\QQ\subset\bar{\bigotimes}_\G\RR=A''$ be given. 
We will follow the proof of Theorem 6 in \cite{solid}.
Passing to a von Neumann subalgebra if necessary, 
we may assume that $\QQ$ is generated by a single unitary 
operator $v\in A''$ with $\lim_k v^k=0$ ultraweakly. 
Fix a non-trivial ultrafilter $\omega$ and consider 
the proper conditional expectation $\Psi_{\QQ}$ from 
$\B(\ell_2\G\otimes L^2A)$ onto $\tilde{\pi}_\tau(\QQ)'$ 
given by 
\[
\Psi_{\QQ}(x)=\mbox{weak$^*$-}\lim_{n\to\omega} 
n^{-1}\sum_{k=1}^n \Ad_{\tilde{\pi}_\tau(v^k)}(x). 
\]
Then, for any 
$x=\sum_{t\in\G}e(t)\otimes x(t)\in \ell_\infty(\G,\K(L^2A))$, 
we have 
\[
\Psi_{\QQ}(x)=\mbox{weak$^*$-}\lim_{n\to\omega} 
\sum_{t\in\G}e(t)\otimes u(t)^*
\big(n^{-1}\sum_{k=1}^n \Ad_{\pi_\tau(v^k)}(u(t)x(t)u(t)^*)\big)u(t)
=0 
\]
since each $u(t)x(t)u(t)^*$ is a compact operator. 
In particular, $\Psi_{\QQ}(K)=0$. 
It follows from Lemma~\ref{lem:inj} and Proposition \ref{prp:kas} 
that the relative commutant $\QQ'\cap\MM$ is injective. 
\end{proof}

\begin{rem}
A similar proof applies to $\MM=\G\ltimes\RR$, where 
$\G\in\CC$ and $\RR$ is the hyperfinite type~$\mathrm{II}_1$ factor. 
(Even though the $C^*$-algebra $C^*_{\red}(\G,\RR)$ is not exact, 
Lemma \ref{lem:inj} is applicable since 
the pair $C^*_{\red}(\G,\RR)\subset\MM$ satisfies local reflexivity.)
It follows that if $\QQ\subset\MM$ is a non-injective subfactor, 
then the relative commutant $\QQ'\cap\MM$ is injective.
In particular, any non-McDuff subfactor 
of $\MM$ is prime. 
This result applies to the factors constructed in \cite{nps}.

Consider an essentially-free measure-preserving ergodic action 
of a non-amenable hyperbolic group $\G$ 
on the standard probability space $[0,1]$ 
and let 
\[
\A=L^\infty[0,1]\subset\G\ltimes L^\infty[0,1]=\MM.
\]
Then, Adams' theorem \cite{adams} says 
$(\A\subset\MM)\neq(\A_1\subset\MM_1)\bar{\otimes}(\A_2\subset\MM_2)$ 
for any type~$\mathrm{II}_1$ factors $\MM_i$ 
with Cartan subalgebras $\A_i\subset\MM_i$. 
Combined with Popa's theorem \cite{popab}, 
this implies that $\MM\neq\MM_1\bar{\otimes}\MM_2$ for any 
HT type~$\mathrm{II}_1$ factors $\MM_i$.
Theorem \ref{thm:crp} generalizes these to 
that $\MM\neq\MM_1\bar{\otimes}\MM_2$ for any 
type~$\mathrm{II}_1$ factors $\MM_i$.

By considering an ergodic but not strongly-ergodic action 
of $\G\in\CC$, we obtain a semisolid (and hence prime) 
type~$\mathrm{II}_1$ factor with the property $(\G)$.

Voiculescu \cite{voiculescuIII} proved that 
the free group factors $L\F_r$ do not have 
regular diffuse injective subalgebras. 
In particular, $L\F_r$ do not have Cartan subalgebras. 
The Bernoulli shift $\F_r\ltimes(\bar{\bigotimes}_{\F_r}\RR)$ 
of $\RR$ by $\F_r$ is an example of a solid factor with 
a Cartan subalgebra. 
\end{rem}

\end{document}